\documentclass[reqno, 11pt]{amsart}
\usepackage{amsfonts,latexsym}
\usepackage{amsmath}
\usepackage{amscd}
\usepackage{float,amsmath,amssymb,mathrsfs,bm,multirow,graphics}
\usepackage[dvips]{graphicx}
\usepackage[percent]{overpic}
\usepackage[all]{xy}

\newcommand{\nequation}{\setcounter{equation}{0}}

\newcommand{\R}{{\Bbb R}}

\newcommand{\C}{{\Bbb C}}

\newcommand{\Z}{{\Bbb Z}}
\newcommand{\CP}{{\Bbb{CP}}}
\newcommand{\proofbegin}{\noindent{\it Proof.\,\,}}
\newcommand{\proofend}{\hfill$\Box$\bigskip}

\newcommand{\projK}{p}
\newcommand{\projV}{q}

\newcommand{\Diff}{\text{\upshape Diff}}

\newcommand{\id}{\text{\upshape id}}

\newcommand{\sect}{\text{\upshape sec}}

\newtheorem{theorem}{Theorem}[section]
\newtheorem{proposition}[theorem]{Proposition}
\newtheorem{corollary}[theorem]{Corollary}
\newtheorem{lemma}[theorem]{Lemma}

\newtheorem{remark}[theorem]{Remark}

\input epsf
\title[Spheres, K\"ahler geometry, and the Hunter-Saxton system]{\sc Spheres, K\"ahler geometry, and \\ the Hunter-Saxton system}
\author{Jonatan Lenells}
\address{Department of Mathematics, Baylor University, One Bear Place \#97328, Waco, TX 76798, USA.}
\email{Jonatan\_Lenells@baylor.edu}
    
\begin{document}

\begin{abstract} 
\noindent
Many important equations of mathematical physics arise geometrically as geodesic equations on Lie groups. In this paper, we study an example of a geodesic equation, the two-component Hunter-Saxton (2HS) system, that displays a number of unique geometric features.
We show that 2HS describes the geodesic flow on a manifold which is isometric to a subset of a sphere. Since the geodesics on a sphere are simply the great circles, this immediately yields explicit formulas for the solutions of 2HS. 
We also show that when restricted to functions of zero mean, 2HS reduces to the geodesic equation on an infinite-dimensional  manifold which admits a K\"ahler structure. We demonstrate that this manifold is in fact isometric to a subset of complex projective space, and that the above constructions provide an example of an infinite-dimensional Hopf fibration. 
\end{abstract}

\maketitle
\noindent
{\small{\sc AMS Subject Classification (2010)}: 53C21, 58D05, 35Q53.}

\noindent
{\small{\sc Keywords}: Diffeomorphism groups, K\"ahler geometry, curvature, nonlinear PDEs.}


\section{Introduction}\nequation
Several well-known equations of mathematical physics arise geometrically as geodesic equations on Lie groups. The classical example is the motion of a rigid body rotating around its center of gravity: the motion is described by the classical Euler equation, which is the geodesic equation on $SO(n)$ endowed with a left-invariant metric defined by the kinetic energy of the body. Another fundamental example is the Euler equation of ideal hydrodynamics: the particles of a fluid moving in a compact $n$-dimensional Riemannian manifold $M$ trace out a geodesic curve in the Lie group of volume-preserving diffeomorphisms of $M$ equipped with a right-invariant metric defined by the kinetic energy of the fluid \cite{A1966, EM1970}. Several other physically relevant equations admit similar geometric formulations, see \cite{V2008}. Once it has been established that an equation admits a formulation of this type, it is tempting to use geometric intuition in order to better understand the behavior of its solutions; for example, directions of positive or negative curvature are expected to correspond to the existence of stable or unstable perturbations of the motion, respectively.

Arnold demonstrated that the curvature of the group of volume-preserving diffeomorphisms associated with hydrodynamics is negative in some directions while it is positive in others \cite{A1966}. For most other geodesic equations, such as the Korteweg-de Vries equation and Camassa-Holm equations, similar results apply---the curvature is sometimes negative and sometimes positive \cite{M1997, M1998}. 
On the other hand, there are a few notable exceptions for which the curvature is of a definite sign. In particular, the groups associated with the Hunter-Saxton (HS) equation, the integrable equation proposed in \cite{KMLP2011}, as well as the two-component versions of both of these equations all have positive and constant curvature \cite{KMLP2011, K2011, Lsphere, LY2011}. For the first two of these equations, this  property has been `explained' as being a consequence of the fact that the underlying spaces are isometric to subsets of $L^2$-spheres \cite{KMLP2011, Lsphere}. Thus, the equations are just the equations for geodesic flow on a sphere in disguise. Since the geodesics on a sphere are simply the great circles, this immediately yields explicit formulas for the solutions of the equations. 
However, for the two-component versions of these equations, a similar geometric interpretation has, until now, been lacking. 

The purpose of the present paper is to `explain' the constant curvature of the two-component Hunter-Saxton (2HS) system by showing that the underlying space is isometric to (part of) the unit sphere in $L^2(S^1; \C)$. Thus, the geometric picture valid for the Hunter-Saxton equation extends also to its two-component version. In the case of HS, the sphere is the unit sphere inÊ $L^2(S^1; \R)$, while in the case of 2HS it is the unit sphere in $L^2(S^1; \C)$. Since the second component of 2HS is related to the complex phase of the functions in $L^2(S^1; \C)$, the geometric picture associated with 2HS reduces to that of HS when the second component vanishes. Even though we restrict our attention to 2HS in this paper, we expect that the two-component equation proposed in \cite{LY2011} admits a similar geometric interpretation. 

We will also consider the restriction of 2HS to solutions $(u, \rho)$, where $\rho$ has zero mean. After showing that the underlying space in this case is an infinite-dimensional K\"ahler manifold with positive (but non-constant) curvature, we will show that it is in fact isometric to a subset of complex projective space. 

Geometrically, the above constructions provide an example of an infinite-dimensional Hopf fibration. Indeed, the circle Ê$S^1$ acts on  functions in the unit sphere $S^\infty \subset L^2(S^1; \C)$ by multiplication by a constant phase. Since the quotient manifold is the infinite-dimensional complex projective space $\CP^\infty$, it follows that $S^\infty$ fibers over $\CP^\infty$ as follows:
\begin{align}\label{hopffibration}
\xymatrix{S^1 \ar[r] & S^\infty \ar[d] \\
& \CP^\infty }
\end{align}
The 2HS equation is (up to isometry) the geodesic equation on $S^\infty$, while the restriction of 2HS to functions $\rho$ of zero mean is the geodesic equation on $\CP^\infty$; the two equations are related by the Hopf fibration (\ref{hopffibration}).

Let us finally point out that geodesic flows on spheres and Hopf fibrations also arise in the analysis of the classical Kepler problem. (Recall that the Kepler problem consists of determining the motion of two point masses interacting under an inverse square force law.) It was noticed by J. Moser that the flow arising in the Kepler problem restricted to the manifold of constant energy $E < 0$ is equivalent (up to a rescaling of time) to the geodesic flow on a sphere. Moreover, the trajectory space of the covering flow on the universal covering space $S^3$ is a two-dimensional sphere $S^2$, and the corresponding map 
$$S^3 \to S^2 \; : \; q \mapsto \{\text{trajectory of } q\}$$ 
is the classical Hopf fibration \cite{M1970}; see also \cite{M1984}. 

After recalling some preliminaries in Section \ref{prelsec}, we establish that the space associated with 2HS is isometric to part of a sphere in Section \ref{spheresec}. In Section \ref{solutionssec}, we analyze solutions of the initial value problem for 2HS. In Section \ref{globalsec}, global properties of the geodesic flow are investigated.
In Section \ref{kahlersec}, we consider the restriction of 2HS to solutions $(u, \rho)$ where $\rho$ has zero mean and show that it describes geodesic flow on a K\"ahler manifold. In Section \ref{hopfsec}, we show that this K\"ahler manifold is isometric to a subset of complex projective space and that the above constructions provide an example of a Hopf fibration. Section \ref{conclusionsec} contains some concluding remarks.

\section{Preliminaries}\label{prelsec}\nequation
\noindent
{\bf The HS and 2HS equations.}
Let $S^1$ denote the circle of length one. The periodic Hunter-Saxton (HS) equation
\begin{align}\label{HS}
  u_{txx} = - 2u_xu_{xx}-uu_{xxx},\qquad t>0, \  x\in S^1,
\end{align}
arises in the study of nematic liquid crystals, $u(t,x)$ being a real-valued function of a space variable $x$ and a slow time variable $t$ \cite{HS1991}. Geometrically, HS is the equation for geodesic flow on the Lie group $\Diff_0(S^1)$ of diffeomorphisms of the circle $S^1$ with a designated fixed point \cite{KM2003}, endowed with the $\dot{H}_1$ right-invariant metric given at the identity by 
$$\langle u, v \rangle_{\id} = \frac{1}{4}\int_{S^1} u_x v_x dx.$$
The space $\Diff_0(S^1)$ equipped with the $\dot{H}^1$-metric is isometric to a subset of the unit sphere in $L^2(S^1; \R)$ \cite{Lsphere}. 

The two-component Hunter-Saxton (2HS) system
\begin{align}\label{2HS}
  \begin{cases}
  u_{txx} = - 2u_xu_{xx}-uu_{xxx} + \rho\rho_x,
  	\\
  \rho_t = -(\rho u)_x,	\end{cases} \qquad t>0, \  x\in S^1,
\end{align}
where Ê$u(t,x)$ and $\rho(t,x)$ are real-valued functions, is a natural generalization of (\ref{HS}). Just like the Hunter-Saxton equation, 2HS is an integrable system with an associated Lax pair formulation and a bi-Hamiltonian structure cf. \cite{CI2008, W2009, W2011}.
Geometrically, (\ref{2HS}) is the equation for geodesic flow on the semidirect product Lie group $G= \Diff_0(S^1) \circledS \mathcal{F}(S^1; S_{4\pi}^1)$, where $S_{4\pi}^1$ denotes the circle of length $4\pi$, $\mathcal{F}(S^1; S_{4\pi}^1)$ denotes the space of (sufficiently smooth) maps $\alpha:S^1 \to S_{4\pi}^1$, and the group $G$ is endowed with the right-invariant metric given at the identity by
\begin{align}\label{metricatid}  
  \langle (u, \rho), (v, \tau) \rangle_{(\id, 0)} = \frac{1}{4}\int_{S^1}\left(u_x v_x + \rho \tau\right) dx.
\end{align}

\noindent
{\bf Diffeomorphism groups.}
In order to set the stage for the rigorous study of (\ref{2HS}) as a geodesic equation, we need to introduce some notation.
Let $s > 5/2$. Let $\Diff^s(S^1)$ denote the Banach manifold of orientation-preserving diffeomorphisms of $S^1$ of Sobolev class $H^s$. We let $\Diff_0^s(S^1)$ denote the subgroup of $\Diff^s(S^1)$ consisting of diffeomorphisms $\varphi$ that keep the point $0 \in S^1 \simeq [0,1)$ fixed, i.e.
$$\Diff_0^s(S^1) = \left\{\varphi \in \Diff^s(S^1) \; | \; \varphi(0) =0\right\}.$$ 
Let $H^s(S^1; \R)$ and $H^s(S^1; \C)$ denote the Hilbert spaces of real-valued and complex-valued functions on $S^1$ of Sobolev class $H^s$, respectively.
Using the identification
\begin{equation}\label{Mhyperplane}  
  \Diff_0^s(S^1) = \{\id + h \, \bigl| \, h \in H^s(S^1; \R), h_x > -1, h(0) = 0\},
\end{equation} 
we can view $\Diff_0^s(S^1)$ as an open subset of the closed hyperplane $\id + H_0^s(S^1; \R) \subset H^s(S^1; \R)$, where $H_0^s(S^1; \R)$ is the closed linear subspace 
$$H_0^s(S^1; \R) = \{u \in H^s(S^1; \R)\bigl| \, u(0) = 0\}.$$ 
Thus, (\ref{Mhyperplane}) provides a global chart for the Banach manifold $\Diff_0^s(S^1)$.

Let $H^{s-1}(S^1; S_{4\pi}^1)$ consist of all maps $S^1 \to S_{4\pi}^1$ of Sobolev class $H^{s-1}$. Let $G^s$ denote the semidirect product $\Diff_0^s(S^1) \circledS H^{s-1}(S^1; S_{4\pi}^1)$ with multiplication given by 
$$(\varphi, \alpha)(\psi, \beta) = (\varphi\circ\psi, \beta + \alpha \circ \psi),$$
where $\circ$ denotes composition and the addition in the second component is pointwise addition of angles, i.e. the addition takes place in $\R/4\pi\Z \simeq S_{4\pi}^1$. $H^{s-1}(S^1; S_{4\pi}^1)$ is a Banach manifold modeled on the space $H^{s-1}(S^1; \R)$; it is the disjoint union of a countable number of components distinguished by the winding number of their elements. It follows that $G^s$ also is a Banach manifold.
The neutral element of $G^s$ is
$(\id,0)$ and $(\varphi,\alpha)$ has the inverse
$(\varphi^{-1},-\alpha \circ\varphi^{-1})$.
The metric $\langle \cdot, \cdot \rangle$ on $G^s$ is defined at the identity by (\ref{metricatid}) and extended to all of $G^s$ by right invariance, i.e.
\begin{align}\label{metric}
\langle U, V \rangle_{(\varphi,\alpha)} & = \left\langle (U_1 \circ\varphi^{-1}, U_2 \circ\varphi^{-1}), (V_1 \circ\varphi^{-1}, V_2 \circ\varphi^{-1}) \right\rangle_{(\text{id},0)}
	\\ \nonumber
& = \frac{1}{4} \int_{S^1} \left(\frac{U_{1x}V_{1x}}{\varphi_x} + U_2V_2 \varphi_x \right) dx,
\end{align}
where $U = (U_1, U_2)$ and Ê$V = (V_1, V_2)$ are elements of $T_{(\varphi, \alpha)}G^s \simeq H_0^{s}(S^1; \R) \times H^{s-1}(S^1; \R)$. 

\begin{remark}\upshape
  1. 
  Equation (\ref{2HS}) is the geodesic equation on $G^s$ in the sense that a curve $(\varphi(t), \alpha(t))$ in $G^s$ is a geodesic if and only if $(u(t), \rho(t)) \in T_{(\id, 0)}G^s$ defined by
\begin{align}\label{urhovarphitft}
 (u,\rho) 
 = (\varphi_t\circ\varphi^{-1},\alpha_t\circ\varphi^{-1})
\end{align}
satisfies (\ref{2HS}), see Proposition \ref{geodesicprop} below.

2. The Hunter-Saxton equation (\ref{HS}) was first shown to arise as a geodesic equation on the homogeneous space $\Diff(S^1)/S^1$ in \cite{KM2003}. Here, we choose to work on the manifold $\Diff_0^s(S^1)$ which is diffeomorphic to $\Diff(S^1)/S^1$; this way we avoid dealing with the coset structure of $\Diff(S^1)/S^1$. 
  
3. As far as regularity is concerned, the geometry of equations (\ref{HS}) and (\ref{2HS}) can be developed in a number of different settings. Here we have chosen to work in the category of Sobolev spaces. Other possible choices for $G$ are
\begin{align*}
  C^nG := C_0^{n}\Diff(S^1) \circledS C^{n-1}(S^1; S_{4\pi}^1)
\end{align*}  
and 
\begin{align*} 
  C^\infty G := C_0^\infty\Diff(S^1) \circledS C^\infty(S^1;S_{4\pi}^1),
\end{align*}
where $C_0^{n}\Diff(S^1)$, $n \geq 2$, denotes the space of orientation-preserving diffeomorphisms of $S^1$ of class $C^n$ that fix $0 \in S^1$. Note that $G^s$ and $C^nG$ are Banach manifolds, but not Lie groups, since left multiplication is only continuous and not smooth. On the other hand, $C^\infty G$ is a Lie group, but only a Fr\'echet manifold. 

4. $G^s$ is a {\it weak} Riemannian manifold in the sense that the topology defined by the metric on each tangent space is weaker than the topology defined by the manifold structure. 

5. It is also possible to view 2HS as the geodesic equation on $\Diff_0^s(S^1) \circledS H^{s-1}(S^1; \R)$. However, since $\rho/2$ is naturally interpreted as the tangent of an angle, it seems more natural to work with $G^s$. In any case, these two spaces are  closely related: $\Diff_0^s(S^1) \circledS H^{s-1}(S^1;\R)$ is a $\Z$-sheeted covering of the component of $G^s$ containing $(\id,0)$; the Lie algebras and the Euler equations for geodesic flow of the two spaces are the same; and the group structure of $\Diff_0^s(S^1) \circledS H^{s-1}(S^1; \R)$ descends to that of $G^s$.
\end{remark}

\section{A sphere}\label{spheresec}\nequation
We will prove that the weak Riemannian manifold $(G^s, \langle \cdot, \cdot \rangle)$ is isometric to a subset of the unit sphere in $L^2(S^1; \C)$. Let $S^\infty$ denote the unit sphere in $L^2(S^1; \C)$ and let $S^{\infty, s}$ denote the elements in $S^\infty$ that are of Sobolev class $H^s$, that is,
$$S^{\infty, s} = \left\{f \in H^{s}(S^1; \C)\; \middle|\; \int_{S^1} |f(x)|^2 dx = 1\right\}.$$
 $S^{\infty, s}$ is a Banach manifold modeled on the closed subspace $1^\perp \subset H^s(S^1; \C)$ of functions orthogonal to the constant function $1 \in H^s(S^1; \C)$,
$$1^\perp := \left\{h \in H^s(S^1; \C) \;\middle |\;  \langle h, 1\rangle_{L^2} = 0 \right\},$$
where $\langle h, 1\rangle_{L^2} = \text{Re} \int_{S^1} h(x) dx$ denotes the component of $h$ along the space spanned by $1$.
Indeed, let
\begin{align}\label{sSdef}
 \sigma_S : S^{\infty, s}\setminus \{-1\} \to 1^\perp, \qquad
\sigma_S(f) = \frac{f - \langle f, 1\rangle_{L^2}}{1 + \langle f, 1\rangle_{L^2}},
\end{align}
denote the stereographic projection from the `south pole' $-1$ with inverse
$$\sigma_S^{-1}(h) = \frac{2h - \|h\|_{L^2}^2 + 1}{\|h\|_{L^2}^2 + 1}.$$
Similarly, define the stereographic projection $\sigma_N$ from the `north pole' $1$ by
$$ \sigma_N : S^{\infty, s}\setminus \{1\} \to 1^\perp, \qquad
\sigma_N(f) = \frac{f - \langle f, 1\rangle_{L^2}}{1 - \langle f, 1\rangle_{L^2}}.$$
Together the two charts defined by $\sigma_S$ and $\sigma_N$ cover $S^{\infty, s}$ and determine its manifold structure.

Let $\mathcal{U}^{s} \subset L^2(S^1; \C)$ denote the open subset of $S^{\infty, s}$ of nowhere vanishing functions:
\begin{align}\label{Usdef}
\mathcal{U}^{s} = \left\{f \in S^{\infty, s} \;\middle|\; |f(x)| > 0 \; \hbox{for} \; x \in S^1\right\}.
\end{align}
We equip $\mathcal{U}^{s}$ with the manifold structure inherited from $S^{\infty, s}$ and the weak Riemannian metric $\langle \cdot, \cdot \rangle_{L^2}$ inherited from $L^2(S^1; \C)$, i.e.
$$\langle X, Y \rangle_{L^2} = \text{Re} \int_{S^1} X(x) \overline{Y(x)} dx,$$
whenever $X, Y \in T_f\mathcal{U}^{s} \subset L^2(S^1; \C)$.

\begin{theorem}\label{sphereth}
The space $(G^s, \langle \cdot, \cdot \rangle)$ is isometric to a subset of the unit sphere in $L^2(S^1; \C)$. 
 More precisely, for any $s > 5/2$, the map $\Phi:G^s \to \mathcal{U}^{s-1} \subset S^\infty$ defined by
  $$\Phi(\varphi, \alpha) = \sqrt{\varphi_x}e^{i\alpha/2}$$
  is a diffeomorphism and an isometry.
\end{theorem}
\proofbegin
If $f \in \mathcal{U}^{s-1}$, then the function $\varphi(x) = \int_0^x |f(y)|^2 dy$ satisfies
$\varphi(0) = 0$, $\varphi(1) = 1$, $\varphi_x = |f|^2 \in H^{s-1}(S^1; \R)$, and $\varphi_x > 0$, while the function $\alpha(x) = 2\arg{f(x)}$ belongs to $H^{s-1}(S^1; S_{4\pi}^1)$.
Thus, the inverse of $\Phi$ is given explicitly by
\begin{align}\label{Phiinverse}
  \Phi^{-1}(f) = \left(\int_0^x |f(y)|^2 dy, 2\arg{f(x)}\right), \qquad f \in \mathcal{U}^{s-1}.
\end{align}  
This shows that $\Phi$ is bijective. Since both $\Phi$ and $\Phi^{-1}$ are smooth, $\Phi$ is a diffeomorphism.
Using that
$$T_{(\varphi,\alpha)} \Phi(U_1, U_2) =  \frac{1}{2\sqrt{\varphi_x} }\left(U_{1x} + iU_2\varphi_x\right)e^{i\alpha/2},$$
we find that
\begin{align*}
\langle T_{(\varphi,\alpha)} \Phi(U_1, U_2), T_{(\varphi,\alpha)} \Phi(V_1, V_2) \rangle_{L^2(S^1; \C)}
& = \text{Re} \int_{S^1} \frac{1}{4\varphi_x}\left(U_{1x} + iU_2\varphi_x\right)\left(V_{1x} - iV_2\varphi_x\right)dx
	\\
& = \frac{1}{4}\int_{S^1} \left(\frac{U_{1x}V_{1x}}{\varphi_x} + U_2V_2 \varphi_x\right) dx
	\\
& = \langle (U_1, U_2), (V_1, V_2) \rangle_{(\varphi, \alpha)},
\end{align*}
whenever $(U_1, U_2)$ and $(V_1, V_2)$ belong to $T_{(\varphi, \alpha)}G^s$.
This shows that $\Phi$ is an isometry. 

\proofend

It follows immediately from Theorem \ref{sphereth} that the sectional curvature of $G^s$ is constant and equal to one. This result was already proved in a different way in \cite{K2011}.

\begin{corollary}\label{curvcor}
  The space $(G^s, \langle \cdot, \cdot \rangle)$ has constant sectional curvature equal to $1$.
\end{corollary}
\proofbegin
In view of Theorem \ref{sphereth}, it is enough to prove that the unit sphere in $L^2(S^1; \C)$ has constant sectional curvature equal to $1$. As in the finite-dimensional case, this can be proved using the Gauss-Codazzi formula. Indeed, letting $n$ denote the outward normal to the sphere, the second fundamental form $\Pi$ is given by
$$\Pi(X, Y) = - \langle \nabla_X n, Y \rangle_{L^2} n = - \langle X, Y \rangle_{L^2} n,$$
where $X,Y$ are vector fields on $S^\infty$. Consequently, if $X$ and $Y$ are orthonormal, the curvature tensor $R$ on the unit sphere satisfies
\begin{align*}
\langle R(X,Y)Y, X \rangle_{L^2}
& = \langle \Pi(X,X), \Pi(Y,Y) \rangle_{L^2}
- \langle \Pi(X,Y), \Pi(Y,X) \rangle_{L^2}.
	\\
& = \langle X,X \rangle_{L^2} \langle Y,Y \rangle_{L^2}
- \langle X, Y \rangle_{L^2} \langle Y,X \rangle_{L^2}  = 1.	
\end{align*}
\proofend

By pulling back the covariant derivative on the sphere $S^\infty$, we can determine the metric connection on $G^s$.
Let $A = -\partial_x^2$. Then $A$ is an isomorphism 
$$H_0^s(S^1; \R) \to \left\{f \in H^{s-2}(S^1; \R) \middle| \; \int_{S^1} f dx = 0\right\}.$$ 
Let $A^{-1}$ be its inverse given by
$$(A^{-1}f)(x) = -\int_0^x\int_0^y f(z)dzdy +  x \int_{S^1}\int_0^y f(z)dz dy$$ 
whenever $\int_{S^1} f dx = 0$. 

\begin{corollary}\label{covdercor}
  The metric covariant derivative on $G^s$ is given by
\begin{align}\label{nabladef}
(\nabla_X Y)(\varphi, \alpha) = DY(\varphi, \alpha) \cdot X(\varphi, \alpha) - \Gamma_{(\varphi, \alpha)}(Y(\varphi, \alpha), X(\varphi, \alpha)),
\end{align}
where the Christoffel map $\Gamma$ is defined for $u = (u_1, u_2)$, $v = (v_1, v_2) \in T_{(\id, 0)}G^s$ by
\begin{subequations}\label{Gammadef}
\begin{align}
\Gamma_{(\id,0)}(u, v) = - \frac{1}{2}\begin{pmatrix} A^{-1}\partial_x(u_{1x}v_{1x} + u_2v_2)) \\
u_{1x}v_2 + v_{1x}u_2 \end{pmatrix},
\end{align}
and extended to all of $G^s$ by right invariance:
\begin{align}
\Gamma_{(\varphi,\alpha)}(u\circ \varphi, v \circ \varphi) = \Gamma_{(\id, 0)}(u, v) \circ \varphi.
\end{align}
\end{subequations}
\end{corollary}
\proofbegin
Let $\nabla'$ be the metric connection on $S^{\infty, s}$. The metric covariant derivative $\nabla$ on $G^s$ is the pull-back of $\nabla'$ by $\Phi$, i.e.
$$\nabla_X Y= \Phi^*(\nabla'_{\Phi_*X}\Phi_*Y).$$
Right invariance implies that it is enough to verify (\ref{nabladef}) at the identity $(\id, 0)$.
We have
$$(\Phi_*Y)(f) = \frac{1}{2}\left(\frac{(Y_1(\Phi^{-1}(f)))_x}{\bar{f}} + i Y_2(\Phi^{-1}(f)) f\right) $$
and
$$(\nabla'_Z W)_f = (DW(f) \cdot Z(f))^t$$
where $Z \mapsto Z^t = Z - \langle Z, f \rangle_{L^2} f$ is the orthogonal projection of $Z$ onto $T_fS^{\infty, s}$.
Thus,
\begin{align}\label{nablaprimePhistar}
(\nabla'_{\Phi_*X}\Phi_*Y)_f
= \frac{1}{2}\left\{\frac{(DY_1 \cdot X)_x}{\bar{f}} + i (DY_2 \cdot X) f
-\frac{Y_{1x}\overline{\Phi_*X}}{\bar{f}^2} + iY_2 \Phi_*X \right\}^t.
\end{align}
Let $u=(u_1, u_2)$ and $v=(v_1, v_2)$ be the values of the vector fields $X$ and $Y$ at the identity, respectively. 
Evaluation of (\ref{nablaprimePhistar}) at $f =1$ yields
\begin{align*}
(\nabla'_{\Phi_*X}\Phi_*Y)_1
= \;& \frac{1}{2}\left\{(DY_1 \cdot u)_x + iDY_2 \cdot u - \frac{1}{2}v_{1x}(u_{1x} - iu_2) + \frac{i}{2}v_2 (u_{1x} + iu_2) \right\}^t
	\\
 =\;& \frac{1}{2}(DY_1 \cdot u)_x  - \frac{1}{4}( u_{1x}v_{1x} + u_2v_2) + \frac{1}{4}\int_{S^1}( u_{1x}v_{1x} + u_2v_2) dx
	\\
& + \frac{i}{2}DY_2 \cdot u + \frac{i}{4}(u_{1x}v_2 + v_{1x}u_2).
\end{align*}
It follows that
\begin{align*}
(\nabla_X Y)(\id, 0) & = \begin{pmatrix} DY_1 \cdot u  - \frac{1}{2}\int_0^x ( u_{1x}v_{1x} + u_2v_2)dx + \frac{x}{2}\int_{S^1}( u_{1x}v_{1x} + u_2v_2) dx  \\
DY_2 \cdot u + \frac{1}{2}(u_{1x}v_2 + v_{1x}u_2)
\end{pmatrix}
	\\
& = DY \cdot u  - \Gamma_{(\id,0)}(u, v),
\end{align*}
which proves (\ref{nabladef}).
\proofend

\begin{remark}\upshape\label{nablasprayremark}
1. It can be verified by direct computation that the $\nabla$ defined in (\ref{nabladef}) defines a covariant derivative on $G^s$ which is compatible with the metric in the sense that
$$Z \langle X, Y \rangle = \langle \nabla_Z X, Y \rangle + \langle X, \nabla_Z Y \rangle$$
for all vector fields $X,Y,Z$ on $G^s$. This gives an alternative proof of Corollary \ref{covdercor}.

2. 
The Christoffel map  (\ref{Gammadef}) defines a smooth spray on $G^s$, i.e., the map
\begin{align*}
(\varphi,\alpha)\mapsto \Gamma_{(\varphi,\alpha)}: G^s \to L^2_{\text{\rm sym}}\left(H_0^{s}(S^1;\R)\times H^{s-1}(S^1; \R); H_0^{s}(S^1;\R)\times H^{s-1}(S^1; \R)\right)
\end{align*}
is smooth (see \cite{EKL2011} for a proof in a similar situation). 
Here, if $X$ and $Y$ are Banach spaces, $L^2_{\text{\rm sym}}(X; Y)$ denotes the Banach space of symmetric bilinear maps from $X$ to $Y$.
\end{remark}

By definition, the geodesics on $G^s$ are the solutions $(\varphi(t), \alpha(t))$ of the equation $\nabla_{(\varphi_t, \alpha_t)}(\varphi_t, \alpha_t) = 0,$ i.e.
\begin{align}\label{geodesiceq}
  (\varphi_{tt}, \alpha_{tt}) = \Gamma_{(\varphi, \alpha)}((\varphi_t, \alpha_t), (\varphi_t, \alpha_t)).
\end{align}
Theorem \ref{sphereth} immediately leads to explicit formulas for the geodesics in $G^s$.

\begin{corollary}\label{geodesiccor}
Let $s > 5/2$. Let $(\varphi(t), \alpha(t)) \in C^\infty([0, T_s); G^s)$ be the geodesic in $G^s$ such that $(\varphi(0), \alpha(0)) = (\id, 0)$ and $(\varphi_t(0), \alpha_t(0)) = (u_0, \rho_0) \in T_{(\id, 0)}G^s$ with maximal time of existence $T_s >0$. Then $(\varphi(t), \alpha(t))$ is given by
\begin{align}\label{geodesicsolution}
  (\varphi(t), \alpha(t)) = \Phi^{-1}\left( \cos{ct} + \frac{u_{0x} + i\rho_0}{2c}\sin{ct}\right), \qquad t \in [0, T_s),
\end{align}  
that is,
\begin{subequations}\label{geodesicsolution2}
\begin{align} 
&\varphi(t,x) = \int_0^x \left\{\left(\cos{ct} + \frac{u_{0x}(y)}{2c}\sin{ct}\right)^2 + \left(\frac{\rho_0(y)}{2c}\sin{ct}\right)^2\right\} dy,   
	\\ \label{geodesicsolution2b}
& \alpha(t,x) = 2\arg\left(\cos{ct} + \frac{u_{0x}(x) + i\rho_0(x)}{2c}\sin{ct}\right),
\end{align}
\end{subequations}
where the speed $c > 0$ of the geodesic is given by $c^2 = \frac{1}{4}\int_{S^1} (u_{0x}^2 + \rho_0^2)dx.$
The maximal existence time $T_s$ is independent of $s > 5/2$ in the sense that if $(u_0, \rho_0) \in T_{(\id, 0)}G^r$ with $r > 5/2$, then $T_r = T_s$ for all $s \in (5/2, r)$.
Moreover, the geodesic $(\varphi(t), \alpha(t))$ exists globally (i.e. $T_s = \infty$) if and only if $\rho_0(x) \neq 0$ for all $x \in S^1$.
\end{corollary}
\proofbegin
The geodesic $f(t)$ on the sphere $S^{\infty, s} \subset L^2(S^1; \C)$ starting at the constant function $1$ with initial velocity $f_t(0)$ is the great circle given explicitly by
\begin{align}\label{greatcircle}
f(t) = \cos{c t} + \frac{f_t(0)}{c} \sin{ ct},
\end{align}
where $c = \|f_t(0)\|_{L^2}$ denotes its speed. 
Indeed, viewing $f(t)$ as a curve in $L^2(S^1; \C)$, we have $f_{tt} = - c^2 f$. Hence, the orthogonal projection of $f_{tt}(t)$ onto the tangent space $T_{f(t)}S^{\infty, s} \subset L^2(S^1; \C)$ vanishes for every $t$. By definition of the induced connectionÊ $\nabla'$ on $S^{\infty, s}$, this shows that $\nabla'_{f_t}f_t \equiv 0$.
Equation (\ref{geodesicsolution}) now follows from Theorem \ref{sphereth} and (\ref{greatcircle}). Equation (\ref{geodesicsolution2}) then follows from (\ref{Phiinverse}).

The geodesic $(\varphi(t), \alpha(t))$ persists as long as 
\begin{align}\label{f1cosct}  
  f(t) := \Phi((\varphi(t), \alpha(t)))= \cos{c t} + \frac{u_{0x} + i\rho_0}{2c}\sin{ct}
\end{align} 
remains in the domain $\mathcal{U}^{s-1}$. The maximal existence time $T_s$ is therefore determined by the time at which $f(t)$ hits the boundary of $\mathcal{U}^{s-1}$, i.e.
$$T_s = \inf \left\{ t > 0 \; \middle|\; f(t,x) = 0 \; \text{for some} \;x \in S^1\right\}.$$
It is clear from this expression that $T_s$ is independent of $s > 5/2$. 

In order to characterize the globally defined geodesics, we need to show that $f(t,x) \neq 0$ for all $x \in S^1$ and $t \geq 0$ if and only if $\rho_0(x) \neq 0$ for all $x \in S^1$. Fix $x \in S^1$. Clearly, by (\ref{f1cosct}), $f(t,x) \neq 0$ for all $t$ if $\rho_0(x) \neq 0$. Conversely, if $\rho_0(x)=0$, then
$$f(t,x) = \cos{c t} + \frac{u_{0x}(x)}{2c}\sin{ct},$$
and for any real number $u_{0x}(x)$ there always exists a $t \geq 0$ such that this expression vanishes (take $t = \pi/2c$ if $u_{0x}(x) = 0$ and $t = \frac{1}{c}[2\pi + \arctan(-2c/u_{0x}(x))]$ if $u_{0x}(x) \neq 0$).
\proofend

\section{Solutions of the Hunter-Saxton system}\label{solutionssec}\nequation
The geometric picture developed above yields explicit expressions for the solutions of the Hunter-Saxton system (\ref{2HS}). It turns out that there exist solutions of 2HS that break in finite time as well as solutions that exist globally. More precisely, we will show that a solution with initial data $(u_0, \rho_0)$ breaks in finite time if and only if $\rho_0(x)$ vanishes at some $x \in S^1$. 

We can write 2HS in the following form suitable for the formulation of weak solutions:
\begin{align}\label{weak2HS}
  & \begin{pmatrix} u_t + uu_x \\ \rho_t + u\rho_x \end{pmatrix} 
  = \begin{pmatrix} 
  - \frac{1}{2} A^{-1}\partial_x\bigl(u_x^2 + \rho^2\bigr) \\ -\rho u_x
  \end{pmatrix}.
\end{align}

\begin{proposition}\label{geodesicprop}
Let $s > 5/2$. Let $(\varphi, \alpha):J \to G^s$ be a $C^2$-curve where $J \subset \R$ is an open interval and define $(u, \rho)$ by (\ref{urhovarphitft}). Then 
\begin{align}\label{urhoCC1J}
(u, \rho) \in C\left(J; H_0^s(S^1;\R) \times H^{s-1}(S^1;\R)\right) \cap C^1\left(J; H_0^{s-1}(S^1;\R) \times H^{s-2}(S^1;\R)\right)
\end{align}
and $(\varphi, \alpha)$ is a geodesic on $J$ if and only if $(u, \rho)$ satisfies the 2HS equation (\ref{weak2HS}) for $t \in J$.
\end{proposition}
\proofbegin
Equation (\ref{urhoCC1J}) follows since, if $q > 3/2$, the composition map $(f, \psi) \mapsto f \circ \psi$ is $C^r$ as a map $H^{q+r}(S^1; \R) \times \Diff^q(S^1) \to H^q(S^1; \R)$, while the inversion map $\psi \mapsto \psi^{-1}$ is $C^r$ as a map $\Diff^{q+r}(S^1) \to \Diff^q(S^1)$ cf. \cite{EM1970}.

Using the right invariance of $\Gamma$, the geodesic equation (\ref{geodesiceq}) can be rewritten as
$$ \begin{pmatrix} u_t + uu_x \\ \rho_t + u\rho_x \end{pmatrix} = \Gamma_{(\id, 0)}((u, \rho), (u, \rho)),$$
which is exactly equation (\ref{weak2HS}).
\proofend

\begin{remark}\upshape
Alternatively, we can use the isometry $\Phi$ of Theorem \ref{sphereth} to show that 2HS is the geodesic equation on $G^s$. Indeed, if $(\varphi(t), \alpha(t))$ is a curve in $G^s$, then $f(t) := \Phi(\varphi(t), \alpha(t))$ is a curve in $\mathcal{U}^{s} \subset S^\infty$, which we know is a geodesic iff $f_{tt}(t) = - c^2 f(t)$, where $c = \|f_t\|_{L^2}$ is the constant speed of the geodesic.
Using the formulas
$$\varphi_t = u \circ \varphi, \qquad \alpha_t = \rho \circ \varphi, \qquad \varphi_{tx} = u_x \circ \varphi \varphi_x$$
we deduce that $f(t)$ satisfies
\begin{align*} \nonumber
f(t) &=  \sqrt{\varphi_x}e^{i\alpha/2},  
	\\ \nonumber
f_t(t)& = \frac{1}{2\sqrt{\varphi_x}} (\varphi_{tx} + i\alpha_t\varphi_x)e^{i\alpha/2} = \frac{\sqrt{\varphi_x}}{2} (u_x + i\rho) \circ\varphi e^{i\alpha/2},
	\\ 
f_{tt}(t) 
 & = \frac{u_x \circ \varphi \sqrt{\varphi_x}}{4}(u_x + i\rho) \circ\varphi e^{\frac{i\alpha}{2}}
+\frac{\sqrt{\varphi_x}}{2}\left[ (u_{x} + i\rho)_t \circ\varphi + (u_{x} + i\rho)_x \circ\varphi \varphi_t\right]e^{\frac{i\alpha}{2}}
	\\ \nonumber
& \quad + \frac{\sqrt{\varphi_x}}{2}(u_x + i\rho) \circ\varphi \frac{i f_t}{2} e^{\frac{i\alpha}{2}}
	\\ 
& = \frac{\sqrt{\varphi_x}}{4}\left\{u_x^2 + iu_x\rho + 2 u_{xt} + 2i\rho_t + 2uu_{xx} + 2iu\rho_x + iu_x \rho - \rho^2 \right\} \circ\varphi  e^{i\alpha/2}.
\end{align*}
It follows that $f_{tt} = - c^2 f$ if and only if
\begin{align}\label{2HSintegrated}
  \begin{cases}
  u_{tx} = - uu_{xx}- \frac{1}{2}u_x^2 + \frac{1}{2}\rho^2 - 2c^2,
  	\\
  \rho_t = -(\rho u)_x,	\end{cases}
\end{align}
where the speed $c$ is given by
$$c^2 = \|(u, \rho)\|_{(\id, 0)}^2 = \frac{1}{4}\int_{S^1} (u_{0x}^2 + \rho_0^2)dx.$$
Since (\ref{2HSintegrated}) is the integrated version of 2HS, this confirms that 2HS is the geodesic equation on $G^s$.
\end{remark}

\begin{theorem}\label{2HSth}
Let $s > 9/2$. Let $(u_0, \rho_0) \in H_0^s(S^1;\R) \times H^{s-1}(S^1;\R)$. Then there exists a unique solution $(u, \rho)$ of the 2HS equation (\ref{weak2HS}) such that 
\begin{align}\label{urhoCC1}
&(u, \rho) \in C\left([0, T_s); H_0^s(S^1;\R) \times H^{s-1}(S^1;\R)\right) 
	\\ \nonumber
&\hspace{3cm} \cap C^1\left([0, T_s); H_0^{s-1}(S^1;\R) \times H^{s-2}(S^1;\R)\right),
	\\ \nonumber
& (u(0), \rho(0)) = (u_0,\rho_0),
\end{align}
where $T_s >0$ is the maximal existence time. The solution is given by
\begin{align}\label{urhovarphialpha}
(u, \rho) = (\varphi_t \circ \varphi^{-1}, \alpha_t \circ \varphi^{-1}), \qquad t \in [0, T_s),
\end{align}
where the curve $(\varphi(t), \alpha(t))$ in $G^s$ is given explicitly in terms of $(u_0, \rho_0)$ by (\ref{geodesicsolution2}).
The maximal existence time $T_s$ is independent of $s > 7/2$ in the sense that if $(u_0, \rho_0) \in H_0^r(S^1) \times H^{r-1}(S^1)$ with $r > 7/2$, then $T_r = T_s$ for all $s \in (7/2, r)$. 
Moreover, the solution $(u(t), \rho(t))$ exists globally (i.e. $T_s = \infty$) if and only if $\rho_0(x) \neq 0$ for all $x \in S^1$.
\end{theorem}
\proofbegin
It follows from Proposition \ref{geodesicprop} that $(u, \rho)$ as defined in (\ref{urhovarphialpha}) is a solution of 2HS satisfying (\ref{urhoCC1}) which exists at least as long as the geodesic $(\varphi, \alpha)$ does. The theorem will follow from Corollary \ref{geodesiccor} if we can show that the maximal existence time $T_s$ of the solution $(u, \rho)$ in (\ref{urhoCC1}) is in fact equal to the maximal existence time of $(\varphi, \alpha)$.
Suppose that $(u, \rho)$ is a solution of 2HS satisfying (\ref{urhoCC1}). Then the map 
$$F:(t, (\psi, \beta)) \mapsto (u(t) \circ \psi, \rho(t) \circ \beta):[0, T_s) \times G^{s-2} \to G^{s-2}$$
is $C^1$. Thus, there exists a unique solution of the ODE
$$(\varphi_t, \alpha_t) = F(t, (\varphi, \alpha))$$
such that $(\varphi(0), \alpha(0)) = (\id, 0)$ and $(\varphi, \alpha) \in C^1([0, T_s); G^{s-2})$.
Then $F(t, (\varphi, \alpha)) \in C^1([0, T_s); G^{s-2})$ so that in fact $(\varphi, \alpha) \in C^2([0, T_s); G^{s-2})$.
Since $(u, \rho)$ satisfies 2HS, $(\varphi, \alpha)$ is a $C^2$-geodesic. But since the spray is smooth this implies that $(\varphi, \alpha) \in C^\infty([0, T_s); G^{s-2})$, i.e. the maximal existence time of $(\varphi, \alpha)$ as a geodesic in $G^{s-2}$ is at least $T_s$. Since the maximal existence time of $(\varphi, \alpha)$ is indepedent of $s$, this shows that the existence times of $(\varphi, \alpha)$ in $G^s$ and of $(u, \rho)$ in (\ref{urhoCC1}) coincide.
\proofend

\begin{corollary}
All solutions of 2HS are periodic in time with period $2\pi$. If $(u, \rho)$ is a solution with maximal existence time $T>0$, then either $T = \infty$ or $T < \pi$.  
\end{corollary}
\proofbegin
$T$ is the smallest time for which the corresponding geodesic in $\mathcal{U}^{s-1}$Ê hits the boundary of $\mathcal{U}^{s-1}$. Since $\mathcal{U}^{s-1}$ is invariant under the antipodal map $f \mapsto -f$ on $S^\infty$, it follows that this happens for $t < \pi$ or it does not happen at all.\proofend

\section{Global behavior of geodesics}\nequation\label{globalsec}
The last statement of Corollary Ê\ref{geodesiccor} gives a characterization of the geodesics on $G^s$ starting at the identity that exist for all times. In this section, we will elaborate further on the global properties of geodesics on $G^s$. 

We begin by describing the geodesic flow on the sphere $S^{\infty, s}$. We let $\exp_1:T_1S^{\infty,s} \to S^{\infty,s}$ denote the (Riemannian) exponential map on $S^{\infty, s}$ restricted to the tangent space at the constant function $1$. The next lemma expresses the fact that given any point $f \in S^{\infty, s}$, there exists a unique great circle passing through $1$ and $f$; unless $f = \pm 1$ in which case there exists an infinite number of such great circles. Thus, the geodesic flow on $S^{\infty, s}$ behaves as can be expected by analogy with the finite-dimensional case.

\begin{lemma}\label{expspherelemma}
  The exponential map $\exp_1:T_1S^{\infty,s} \to S^{\infty,s}$ on the sphere $S^{\infty, s}$ satisfies
\begin{align*}
\exp_{1}^{-1}(f) = \begin{cases}
 \{(r_0 + 2\pi n)X_0 \; | \; n \in \Z\}, & f \in S^{\infty, s} \setminus \{1, -1\}, \\
\left\{rX \in T_1S^{\infty,s} \; |\; \|X\|_{L^2} = 1, r \in 2\pi \Z\right\}, & f = 1, \\
 \left\{rX \in T_1S^{\infty,s} \; |\; \|X\|_{L^2} = 1, r \in \pi + 2\pi \Z\right\}, & f = -1, 
 \end{cases}
 \end{align*}
where the unit length vector $X_0 \in T_1S^{\infty,s}$ and the real number $r_0\in (0, \pi)$ are given by
\begin{align}\label{X0r0def}
X_0 = \frac{f - \langle f , 1 \rangle_{L^2}}{\sqrt{1 - \langle f, 1\rangle_{L^2}^2}},
\qquad r_0 = \arccos \langle f, 1 \rangle_{L^2}.
\end{align}
\end{lemma}
\proofbegin
Let $f \in S^{\infty,s}$.
If $X \in T_1S^{\infty,s}$ has length one, we have (cf. equation (\ref{greatcircle}))
$$\exp_1(rX) = \cos{r} + X \sin{r}.$$
Thus, $\exp_1(rX) = f$ if and only if
\begin{align}\label{cosrXsinrf}
\cos{r} + X \sin{r} = f.
\end{align}
Applying $\langle \cdot,1 \rangle_{L^2}$ to both sides of this equation, we find
$$\cos{r} = \langle f,1 \rangle_{L^2},$$
and the lemma now follows from (\ref{cosrXsinrf}).
\proofend
  
Given two points $f, g \in S^{\infty, s}$ such that $f \neq \pm g$, Lemma \ref{expspherelemma} implies that there is a unique geodesic of length $r_0 \in (0, \pi)$ joining $f$ to $g$; we call this the {\it short geodesic segment from $f$ to $g$}. There is also a unique geodesic of length $2\pi - r_0 \in (\pi, 2\pi)$ joining $f$ to $g$, which goes around the sphere in the opposite direction; we call this the {\it long geodesic segment from $f$ to $g$}.
If $f, g$ belong to $\mathcal{U}^s$, we may ask whether the short and long geodesic segments connecting $f$ to $g$ are also contained in $\mathcal{U}^s$. Clearly, since $\mathcal{U}^s$ is invariant under the antipodal map $f \mapsto -f$, the short geodesic segment lies in $\mathcal{U}^s$ whenever the long segment does.
  
\begin{proposition}\label{shortlongprop}
Let $f, g \in \mathcal{U}^{s}$ and suppose that $f \neq \pm g$. The short geodesic segment from $f$ to $g$ is contained in $\mathcal{U}^{s}$ if and only if $f(x)/g(x) \notin (-\infty, 0)$ for all $x \in S^1$.
The long geodesic segment from $f$ to $g$ is contained in $\mathcal{U}^{s}$ if and only if $f(x)/g(x) \notin \R$ for all $x \in S^1$.
\end{proposition}
\proofbegin
Let $\Phi((\varphi, \alpha)) = f$ and $\Phi((\psi, \beta)) = g$ be two points in $\mathcal{U}^{s}$. Right invariance implies that there exists a geodesic from $(\varphi, \alpha)$ to $(\psi, \beta)$ in $G^{s+1}$ iff there exists one from $(\id, 0)$ to $(\varphi, \alpha)(\psi, \beta)^{-1} = (\varphi \circ \psi^{-1}, (\alpha - \beta) \circ \psi^{-1})$. Moreover, right translation preserves the length of a geodesic. 
Hence, the short (long) geodesic segment from $f$ to $g$ is contained in $\mathcal{U}^{s}$ iff the short (long) geodesic segment from $\Phi((\id, 0)) = 1$ to
$$\Phi( (\varphi \circ \psi^{-1}, (\alpha - \beta) \circ \psi^{-1})) = \left(\frac{f}{g}\right) \circ \psi^{-1}$$
is contained in $\mathcal{U}^{s}$.
It is therefore enough to prove the proposition in the case when $g = 1$.

Let $f \in \mathcal{U}^s$ with $f \neq \pm 1$. Let $X_0$ and $r_0$ be as in (\ref{X0r0def}). We claim that
\begin{align}\label{exp1rX0}
\exp_1(rX_0) = \cos{r} + X_0 \sin{r} = \frac{\sin(r_0 - r) + f\sin{r}}{\sin{r_0}}
\end{align}
belongs to $\mathcal{U}^s$ for $r \in (0,r_0)$ if and only if $f(x) \notin (-\infty, 0)$ for $x \in S^1$.
Indeed, let us fix $x\in S^1$. Then $\exp_1(rX_0)(x) = 0$ for some $r \in (0,r_0)$ iff $f(x) = -\sin(r_0 -r)/\sin{r}$
for some $r \in (0,r_0)$.
Since the right-hand side of this equation maps the interval $(0,r_0)$ to $(-\infty, 0)$, we see that $\exp_1(rX_0)(x) = 0$ for some $r \in (0,r_0)$ iff $f(x) < 0$. This proves the first half of the proposition.
In order to prove the second half, we need to show that the geodesic $\exp_1(rX_0)$ lies in $\mathcal{U}^s$ for $r \in \R$ if and only if $f(x) \notin \R$ for $x \in S^1$. This can either be deduced from (\ref{exp1rX0}) or be seen as a consequence of the last statement of Corollary Ê\ref{geodesiccor} using that  $\text{Im}(r_0X_0) = \frac{r_0}{\sin{r_0}} \text{Im}\; f$.
\proofend

Note that the antipodal involution $f \mapsto -f$ on $S^{\infty, s}$ corresponds under the isometry $\Phi$ to the involution $(\varphi, \alpha)\mapsto (\varphi, \alpha + 2\pi)$ of $G^s$.
Thus, if we use the isometry $\Phi$ to transfer the result of Proposition \ref{shortlongprop} to $G^s$, we immediately find the following result.
 
\begin{corollary}\label{Gsgeocor}
Let $(\varphi, \alpha), (\psi, \beta) \in G^s$ be distinct points in $G^s$ and suppose that $(\varphi, \alpha) \neq (\psi, \beta + 2\pi)$. There exists a geodesic joining $(\varphi, \alpha)$ to $(\psi, \beta)$ if and only if $e^{i(\alpha(x) - \beta(x))/2} \neq -1$ for all $x \in S^1$. This geodesic is unique and has length less than $\pi$ provided that there exists an $x$ such that $e^{i(\alpha(x) - \beta(x))/2} = 1$. If no such $x$ exists (so that $e^{i(\alpha(x) - \beta(x))/2} \neq \pm 1$ for allÊ $x \in S^1$), then the geodesic is defined on all of $\R$, is periodic with period $2\pi$ with respect to an arc-length parameter, and is unique up to the choice of its direction.
On the other hand, for any $(\varphi, \alpha) \in G^s$, there exists an infinite number of geodesics joining $(\varphi, \alpha)$ to $(\varphi, \alpha + 2\pi)$. All of these geodesics exist globally and are $2\pi$-periodic with respect to an arc-length parameter.
\end{corollary}


Let $\exp_{(\id, 0)}$ denote the (Riemannian) exponential map on $G^s$ restricted to $T_{(\id, 0)}G^s$. 
The domain $D$ of $\exp_{(\id, 0)}$ consists of all $(u_0, \rho_0)$ such that the geodesic starting at $(\id, 0)$ with initial velocity $(u_0, \rho_0)$ exists for a time larger than $1$. 
Using the above results, it is possible to express $\exp_{(\id, 0)}$ and its multivalued inverse explicitly. The following proposition gives the expression for $\exp_{(\id, 0)}^{-1}((\varphi, \alpha))$ in the case that $(\varphi, \alpha) \neq (\id, 0)$ and $(\varphi, \alpha) \neq (\id, 2\pi)$.

\begin{proposition}\label{logprop}
  Let $(\varphi, \alpha) \in G^s$ and suppose that $(\varphi, \alpha) \neq (\id, 0)$ and $(\varphi, \alpha) \neq (\id, 2\pi)$. Then $\exp_{(\id, 0)}^{-1}((\varphi, \alpha))$ is the empty set if $e^{i\alpha(x)/2} = -1$ for some $x \in S^1$. Assuming that $e^{i\alpha(x)/2} \neq -1$ for all $x \in S^1$, we have
 $$\exp_{(\id, 0)}^{-1}((\varphi, \alpha))
 = \begin{cases} r_0(u_0, \rho_0) & \text{if } e^{i\alpha(x)/2} = 1 \text{ for some } x \in S^1, \\
 \left\{(r_0 + 2\pi n)(u_0, \rho_0) \; | \; n \in \Z \right\} &	\text{otherwise},
 \end{cases}$$
where the unit length vector $(u_0, \rho_0) \in T_{(\id,0)} G^s$ and the real number $r_0 \in (0, \pi)$ are given by
\begin{align}\nonumber
& u_0(x) = \frac{2}{\sqrt{ 1- \mu_0^2}} \int_0^x \left(\sqrt{\varphi_x(y)} \cos\bigl(\frac{\alpha(y)}{2}\bigr)- \mu_0\right) dy,
\qquad \rho_0(x) = \frac{2\sqrt{\varphi_x(x)}\sin\bigl(\frac{\alpha(x)}{2}\bigr)}{\sqrt{ 1- \mu_0^2}},
	\\ \label{u0rh0r0}
& r_0 = \arccos \int_{S^1} \sqrt{\varphi_x} \cos\frac{\alpha}{2}dx,
\end{align}
and $\mu_0 = \int_{S^1} \sqrt{\varphi_x} \cos(\alpha/2) dx$.
\end{proposition}
\proofbegin
The expressions in (\ref{u0rh0r0}) follow from Lemma \ref{expspherelemma} since $T_{(\id, 0)}\Phi (u_0, \rho_0) = X_0,$
where $X_0$ is as in (\ref{X0r0def}) with $f = \sqrt{\varphi_x}e^{i\alpha/2}$. The rest follows from Corollary \ref{Gsgeocor}.
\proofend

\begin{remark}\upshape
  If $\alpha \equiv 0$, the above results reduce to those derived in \cite{Lsphere} for HS. Nevertheless, there are big differences between the geometries associated with 2HS and HS. For example, for the Hunter-Saxton equation, any two points of the underlying space can be joined by a unique length-minimizing geodesic \cite{Lsphere}. In contrast, for 2HS we have seen that there are points that can be joined by more than one geodesic as well as points that cannot be joined by any geodesic even though they lie in the same component of $G^s$. 
  
\end{remark}

\section{A K\"ahler manifold}\nequation\label{kahlersec}
The mean value $\int_{S^1} \rho dx$ of the second component $\rho$ of a solution $(u, \rho)$ of 2HS is conserved, i.e.
$$\frac{d}{dt}\int_{S^1} \rho dx = -\int_{S^1} (\rho u)_x dx = 0.$$
Thus, if $\rho$ has zero mean initially, it will have zero mean at all later times. This suggests that we consider the following variation of 2HS:
\begin{align}\label{pi2HS}
  \begin{cases}
  u_{txx} = - 2u_xu_{xx}-uu_{xxx} + \pi(\rho)\rho_x,
  	\\
  \pi(\rho)_t = -(\pi(\rho) u)_x,	\end{cases} \qquad t>0, \  x\in S^1,
\end{align}
where $\pi(\rho) =  \rho - \int_{S^1} \rho dx$ denotes the orthogonal projection onto the subspace of functions of zero mean. 

For solutions such that $\int_{S^1} \rho dx = 0$, equation (\ref{pi2HS}) coincides with 2HS.
However, we will see that (\ref{pi2HS}) possesses some interesting geometric properties not shared by the 2HS equation (\ref{2HS}). In particular, (\ref{pi2HS}) is the geodesic equation on a manifold $K$ which admits a K\"ahler structure. 
In this section, we will introduce the space $K$, show that it is a K\"ahler manifold, and prove that (\ref{pi2HS}) is the associated geodesic equation.

\medskip
\noindent
{\bf The K\"ahler manifold $K^s$.}
Let $s > 5/2$. Let $H^{s-1}(S^1; S_{4\pi}^1)/S_{4\pi}^1$ denote the space $H^{s-1}(S^1; S_{4\pi}^1)$ with two elements being identified iff they differ by a constant phase; the equivalence class of $\alpha \in H^{s-1}(S^1; S_{4\pi}^1)$ will be denoted by $[\alpha] \in H^{s-1}(S^1; S_{4\pi}^1)/S_{4\pi}^1$.
We define $K^s$ as the semidirect product $\Diff_0^s(S^1) \circledS (H^{s-1}(S^1; S_{4\pi}^1)/S_{4\pi}^1)$ with multiplication given by 
$$(\varphi, [\alpha])(\psi, [\beta]) = (\varphi\circ\psi, [\beta + \alpha \circ \psi]).$$
Let $H^{s-1}(S^1; \R)/\R$ denote the space $H^{s-1}(S^1; \R)$ with two functions being identified iff they differ by a constant. Since the constant functions form a closed linear subspace of $H^{s-1}(S^1; \R)$, the quotient space $H^{s-1}(S^1; \R)/\R$ is a Banach space.
The space $H^{s-1}(S^1; S_{4\pi}^1)/S_{4\pi}^1$ is a Banach manifold modeled on $H^{s-1}(S^1; \R)/\R$. Together with the global chart (\ref{Mhyperplane}) for $\Diff_0^s(S^1)$, this turns $K^s$ into a Banach manifold.

We equip $K^s$ with the right-invariant metric given at the identity by
\begin{align}\label{pimetricatid}  
  \langle (u, [\rho]), (v, [\tau]) \rangle_{(\id, [0])} = \frac{1}{4}\int_{S^1} (u_x v_x + \pi(\rho)\pi(\tau)) dx.
\end{align}
Extending the projection $\pi$ to any tangent space by right invariance so that $\pi(U_2) = U_2 - \int_{S^1} U_2 \varphi_x dx$
whenever $(U_1, [U_2]) \in T_{(\varphi, \alpha)}K^s \simeq H_0^s(S^1; \R) \times (H^{s-1}(S^1; \R)/\R)$, we have
\begin{align}\label{pimetric}
\langle (U_1, [U_2]), (V_1, [V_2]) \rangle_{(\varphi, [\alpha])} = \frac{1}{4}\int_{S^1} \left(\frac{U_{1x}V_{1x}}{\varphi_x} + \pi(U_2)\pi(V_2)\varphi_x\right) dx.
\end{align}

We define a connection $\nabla$ on $K^s$ by
$$\nabla_X Y = DY \cdot X - \Gamma_{(\varphi, [\alpha])}(Y, X),$$
where the Christoffel map $\Gamma$ is defined for $u = (u_1, [u_2])$, $v = (v_1, [v_2])$ in $T_{(\id, [0])}K^s$ by
\begin{subequations}\label{piGammadef}
\begin{align}
\Gamma_{(\id, [0])}(u, v) = - \frac{1}{2}\begin{pmatrix} A^{-1}\partial_x(u_{1x}v_{1x} + \pi(u_2)\pi(v_2)) \\
[u_{1x}\pi(v_2) + v_{1x}\pi(u_2)] \end{pmatrix}
\end{align}
and extended to the tangent space at $(\varphi, [\alpha]) \in K^s$ by right invariance:
\begin{align}
\Gamma_{(\varphi, [\alpha])}(u\circ \varphi, v \circ \varphi) = \Gamma_{(\id, [0])}(u, v) \circ \varphi.
\end{align}
\end{subequations}
We also define a (1,1)-tensor $J$ and a two-form $\omega$ on $K^s$ by
\begin{align}\label{Jdef}
J_{(\varphi, [\alpha])}(U_1, [U_2]) = \left(-\int_0^x \pi(U_2)\varphi_x dy, \left[\frac{U_{1x}}{\varphi_x}\right]\right)
\end{align}
and
\begin{align}\label{omegadef}
\omega_{(\varphi, [\alpha])}((U_1, [U_2]), (V_1, [V_2]))
= \frac{1}{4}\int_{S^1}(U_{2x}V_{1} - V_{2x}U_{1}) dx
\end{align}
whenever $(U_1, [U_2]), (V_1, [V_2]) \in T_{(\varphi, [\alpha])}K^s$. 
Note that $\omega$ and $J$ are right-invariant. Indeed, a change of variables in (\ref{omegadef}) shows that
$$\omega_{(\varphi, [\alpha])}(u \circ \varphi, v \circ \varphi) = \omega_{(\id, [0])}(u, v),$$
while right-invariance of $J$ follows by a simple calculation:
\begin{align*}
J \circ TR_{(\varphi, [\alpha])}(u_1, [u_2]) 
& = J(U_1, [U_2]) = \left(-\int_0^x \pi(U_2)\varphi_x dy, \left[\frac{U_{1x}}{\varphi_x}\right]\right)
	\\
& = \left(-\int_0^{\varphi(x)} \pi(u_2) dy, [u_{1x} \circ \varphi]\right)
 = TR_{(\varphi, [\alpha])} \circ J(u_1, u_2),
\end{align*}
if $(U_1, [U_2]) = (u_1, [u_2]) \circ \varphi \in T_{(\varphi, [\alpha])}K^s$ and $R_{(\varphi, [\alpha])}$ denotes right multiplication by $(\varphi, [\alpha])$.
We refer to \cite{Lang} for an introduction to differential forms and tensor fields on Banach manifolds.

\begin{theorem}\label{kahlerth}
$K^s$ is a K\"ahler manifold. In fact, letting $g$ denote the metric $\langle \cdot, \cdot \rangle$ on $K^s$, the following hold:
\begin{itemize}
\item[(a)] $g$ is a smooth metric on $K^s$ and $\nabla$ is a smooth connection compatible with $g$.

\item[(b)] $\omega$ is a symplectic form on $K^s$ compatible with $\nabla$, i.e. $\omega$ is a smooth nondegenerate closed two-form on $K^s$ such that $\nabla \omega = 0$.

\item[(c)] $J$ is a complex structure on $K^s$ compatible with $\nabla$, i.e. $J$ is a smooth $(1,1)$-tensor on $K^s$ such that $J^2 = -I$ and $\nabla J = 0$.

\item[(d)] The symplectic form $\omega$, the metric $g$, and the complex structure $J$ are compatible, i.e. $\omega(U, V) = g(JU, V).$

\item[(e)] The metric $g$ is almost Hermitian, i.e. $g(U, V) = g(JU, JV).$

\item[(f)]  The Nijenhuis tensor $N^J$ defined for vector fields $X, Y$ by
$$N^J(X,Y) = [X,Y] + J[JX, Y] + J[X, JY] - [JX, JY]$$
vanishes identically.
\end{itemize}
\end{theorem}
\proofbegin
Throughout the proof, $u = (u_1, [u_2])$, $v = (v_1, [v_2])$, and $w = (w_1, [w_2])$ will denote elements of $T_{(\id, [0])}K^s$.

\smallskip\noindent
{\bf Proof of (a).}
Smoothness of $g$ follows since (\ref{pimetric}) depends smoothly on $\varphi$ and $\alpha$ as a bilinear map from $H_0^s(S^1; \R) \times (H^{s-1}(S^1; \R)/\R)$ to $\R$.
$\nabla$ defines a smooth connection because, as in the case of $G^s$ above, the Christoffel map $\Gamma$ defined in (\ref{piGammadef}) defines a smooth spray on $K^s$, i.e., the map
\begin{align*}
& (\varphi,[\alpha])\mapsto \Gamma_{(\varphi,[\alpha])}
	\\
& K^s \to L^2_{\text{\rm sym}}\left(H_0^{s}(S^1;\R)\times (H^{s-1}(S^1; \R)/\R);
H_0^{s}(S^1;\R)\times (H^{s-1}(S^1; \R)/\R)\right)
\end{align*}
is smooth, cf. Remark \ref{nablasprayremark}. 

In order to show that $\nabla$ and $g$ are compatible, we need to show that
$$X\langle Y, ZÊ\rangle = \langle \nabla_X Y , Z \rangle + \langle Y, \nabla_X Z \rangle$$
for any vector fields $X,Y,Z$ on $K^s$.
By right invariance, it is enough to verify this identity at the identityÊ $(\id, [0])$.
Moreover, using the argument on pp. 129-130 of \cite{EM1970}, it is enough to prove it when $X, Y, Z$ are right-invariant vector fields. Thus, assume that
\begin{align}\label{XYZuvw}
& X(\varphi, \alpha) = (u_1 \circ \varphi, [u_2 \circÊ\varphi]), \qquad
Y(\varphi, \alpha) = (v_1 \circ \varphi, [v_2 \circÊ\varphi]), 
	\\ \nonumber
& Z(\varphi, \alpha) = (w_1 \circ \varphi, [w_2 \circÊ\varphi]).
\end{align}
Then, the function $\langle Y, ZÊ\rangle $ is constant so that $X\langle Y, ZÊ\rangle = 0$. Moreover,
\begin{align*}
\langle \nabla_X Y , Z \rangle_{(\id, [0])}
 = \;& \langle \left(v_{1x} u_1 - \Gamma_1(v, u), [v_{2x}u_1 - \Gamma_2(v, u)]\right) , (w_1, [w_2]) \rangle_{(\id, 0)}
	\\
 = \;& \frac{1}{4}\int_{S^1} \biggl\{\left(v_{1x}u_1 + \frac{1}{2}A^{-1}\partial_x(u_{1x} v_{1x} + \pi(u_2)\pi(v_2))\right)_x w_{1x}
	\\
& \hspace{1.3cm} + \pi\left(v_{2x}u_1 + \frac{1}{2}(u_{1x}\pi(v_2) + v_{1x}\pi(u_2))\right)\pi(w_2)\biggr\}dx
	\\
 =\;& \frac{1}{4}\int_{S^1} \biggl\{(v_{1x}u_1)_xw_{1x} - \frac{1}{2}(u_{1x} v_{1x} + \pi(u_2)\pi(v_2)) w_{1x} + v_{2x}u_1\pi(w_2) 
 	\\
& \hspace{1.3cm}+ \frac{1}{2}(u_{1x}\pi(v_2) + v_{1x}\pi(u_2))\pi(w_2) \biggr\}dx.
\end{align*}
Simplification using integration by parts leads to 
\begin{align*}
\langle \nabla_X Y , Z \rangle_{(\id, [0])} &+ \langle Y, \nabla_X Z \rangle_{(\id, [0])}
	\\
& = \frac{1}{4} \int_{S^1} \biggl\{  v_{2x}u_1\pi(w_2) + u_{1x}\pi(v_2)\pi(w_2) + w_{2x}u_1\pi(v_2)  \biggr\}dx = 0.
\end{align*}

\smallskip\noindent
{\bf Proof of (b).} First note that $\omega$ is a smooth two-form on $K^s$, because the right-hand side of (\ref{omegadef}) is skew-symmetric in $U, V$ and indepedent of $(\varphi, [\alpha])$. Moreover, $\omega$ is non-degenerate because if $\omega((U_1, [U_2]), (V_1, [V_2])) = 0$ for all $V_1, V_2$, then $U_1 = U_{2x} = 0$, which means that $(U_1, [U_2]) = 0$.
In order to show that $\omega$ is closed (i.e. $d\omega = 0$), we recall that the exterior derivative $d\eta$ of a $p$-form $\eta$ on a Banach manifold is defined by
\begin{align}\label{ddef}
d\eta(X_0, \dots, X_p) = \; & \sum_{i = 0}^p (-1)^i X_i(\eta(X_0, \dots, \hat{X}_i, \dots, X_p))
	\\ \nonumber
& + \sum_{0 \leq i < j \leq p} (-1)^{i+j}\eta([X_i, X_j], X_0, \dots, \hat{X}_i, \dots, \hat{X}_j, \dots, X_p),
\end{align}
so that in a local chart 
$$d\eta(X_0, \dots, X_p) = \sum_{i=0}^p (-1)^i (D\eta \cdot X_i)(X_0, \dots, \hat{X}_i, \dots, X_p).$$
Since the local expression (\ref{omegadef}) of $\omega$Ê is independent of $(\varphi, [\alpha])$, it follows that $d\omega = 0$.

It remains to prove that $\nabla \omega = 0$. We need to verify that
\begin{align}\label{omegaparallel}
  X \omega(Y, Z)  = \omega( \nabla_X Y, Z ) + \omega( Y, \nabla_X Z ),
\end{align}
for all vector fields $X,Y,Z$. As in the proof of the identity $\nabla g = 0$, it is enough to verify (\ref{omegaparallel}) at the identity and in the case that $X,Y,Z$ are right invariant. Thus, let $X,Y,Z$ be right-invariant vector fields as in (\ref{XYZuvw}). Then $X \omega(Y, Z) = 0$ and
\begin{align*}
\omega_e ( \nabla_X Y, Z )
& = \omega( \left(v_{1x} u_1 - \Gamma_1(v, u), [v_{2x}u_1 - \Gamma_2(v, u)]\right), (w_1, [w_2]) )
	\\
& = \frac{1}{4}\int_{S^1}\left\{(v_{2x}u_1 - \Gamma_2(v, u))_x w_{1} - w_{2x}\left(v_{1x} u_1 - \Gamma_1(v, u)\right)\right\}dx
	\\
& = \frac{1}{4}\int_{S^1}\biggl\{\left(v_{2x}u_1 + \frac{1}{2}(u_{1x}\pi(v_2) + v_{1x}\pi(u_2))\right)_x w_1 
	\\
&\hspace{3cm}- \pi(w_2)_x \left(v_{1x} u_1 + \frac{1}{2}A^{-1}\partial_x(u_{1x}v_{1x} + \pi(u_2)\pi(v_2))\right)\biggr\}dx.
\end{align*}
In view of the identity $\partial_x A^{-1}\partial_x f = -f + \int_{S^1} f dx,$
we can write the right-hand side as
\begin{align*}
&\frac{1}{4}\int_{S^1}\biggl\{-u_1(v_{2x}w_{1x} + w_{2x}v_{1x})
  - \frac{1}{2}(u_{1x}\pi(v_2) + v_{1x}\pi(u_2)) w_{1x} 
  	\\
& \hspace{5cm}   - \frac{1}{2} \pi(w_2) (u_{1x}v_{1x} + \pi(u_2)\pi(v_2))\biggr\}dx
\end{align*}
Since the part of this expressions that is antisymmetric in $v$ and $w$ vanishes, we deduce that
$\omega( \nabla_X Y, Z ) + \omega( Y, \nabla_X Z ) = 0.$

\smallskip\noindent
{\bf Proof of (c).}
Smoothness of $J$ follows since the right-hand side of (\ref{Jdef}) depends smoothly on $\varphi$ as a linear map from $H_0^s(S^1; \R) \times (H^{s-1}(S^1; \R)/\R)$ to itself. Since
$$\pi\Bigl(\frac{U_{1x}}{\varphi_x}\Bigr) = \frac{U_{1x}}{\varphi_x},$$
we have $J^2(U_1, [U_2]) = -(U_1, [U_2]),$
showing that $J$ is an almost complex structure.
It only remains to prove that $\nabla J = 0$. This can either be seen as a consequence of $\nabla g = \nabla \omega = 0$ and the statement $(d)$ proved below, or be established directly as follows. The covariant derivative of the $(1,1)$-tensor $J$ is given locally by
$$(\nabla_X J)(Y) = (DJ\cdot X)(Y) - \Gamma(JY, X) + J\Gamma(Y, X).$$
We compute each of the terms in turn:
\begin{align*} \nonumber
 (DJ(\varphi&, \alpha) \cdot X)(Y) 
	\\ \nonumber
& = \frac{d}{d\epsilon}\biggl|_{\epsilon = 0}
\biggl(- \int_0^x \left(Y_2 - \int_{S^1} Y_2 (\varphi_x + \epsilon X_{1x})dz \right) (\varphi_x + \epsilon X_{1x}) dy, 
\left[\frac{Y_{1x}}{\varphi_x + \epsilon X_{1x}}\right]\biggr)
	\\ 
& = \left(- \int_0^x \pi(Y_2) X_{1x} dy + \varphi(x) \int_{S^1}Y_{2}X_{1x}dy , - \left[\frac{Y_{1x} X_{1x}}{\varphi_x^2}\right]\right).
\end{align*}
Since $A^{-1}\partial_x f = -\int_0^x fdy + x \int_{S^1} fdy,$
evaluation at the identity yields
$$(DJ \cdot u)(v)  = \left(A^{-1}\partial_x(\pi(v_2) u_{1x}), - [v_{1x} u_{1x}]\right).$$
where $(u_1, [u_2])$ and $(v_1, [v_2])$ are the values of $X$ andÊ $Y$ at the identity, respectively.
Moreover,
$$- \Gamma(Jv, u)
= \frac{1}{2}\begin{pmatrix} A^{-1}\partial_x(-\pi(v_2) u_{1x} + \pi(v_{1x})\pi(u_2)) \\
[- \pi(v_2)\pi(u_2) + u_{1x}\pi(v_{1x})] \end{pmatrix}
$$
and
$$J\Gamma(v, u)
= -\frac{1}{2}\begin{pmatrix}
A^{-1}\partial_x \pi(u_{1x}\pi(v_2) + v_{1x}\pi(u_2)) dx \\
[\partial_xA^{-1}\partial_x(u_{1x}v_{1x} + \pi(u_2)\pi(v_2))]
\end{pmatrix}.$$
The sum of the preceding three equations vanishes; thus $\nabla J = 0$.

\smallskip\noindent
{\bf Proof of (d).}
We have
\begin{align*}
g(JU, V) & = g\left( \left(-\int_0^x \pi(U_2)\varphi_x dx, \left[\frac{U_{1x}}{\varphi_x}\right] \right), V \right)
	\\
& = \frac{1}{4}\int_{S^1} \left( \frac{-\pi(U_2)\varphi_x V_{1x}}{\varphi_x} + \pi\left(\frac{U_{1x}}{\varphi_x}\right)\pi(V_2) \varphi_x\right) dx = \omega(U, V).
\end{align*}

\smallskip\noindent
{\bf Proof of (e).}
This is a simple calculation:
\begin{align*}
g(JU, JV) = 
g\left(\left(-\int_0^x \pi(U_2)\varphi_x dx, \left[\frac{U_{1x}}{\varphi_x}\right]\right), \left(-\int_0^x \pi(V_2)\varphi_x dx, \left[\frac{V_{1x}}{\varphi_x}\right] \right) \right)
	\\
= \frac{1}{4}\int_{S^1} \left\{ \frac{\pi(U_2)\varphi_x \pi(V_2) \varphi_x}{\varphi_x} + \pi\left(\frac{U_{1x}}{\varphi_x}\right)\pi\left(\frac{V_{1x}}{\varphi_x}\right)\varphi_x \right\} dx = g(U, V).
\end{align*}

\smallskip\noindent
{\bf Proof of (f).}
It is enough to verify this at the identity by right invariance. We let $[\cdot, \cdot]$ be the Lie bracket on $T_{(\id, [0])}K^s$ induced by right-invariant vector fields. Then
\begin{align}\label{Ksbracket}
 [u,v] & = \begin{pmatrix} v_{1x}u_1 - u_{1x}v_1 \\ [v_{2x}u_1 - u_{2x}v_1] \end{pmatrix},
	\\\nonumber
J[Ju,v] & = \begin{pmatrix} \int_0^x \pi\left(v_{2x}\int_0^y \pi(u_2)dy + u_{1xx}v_1\right) dy \\
\left[(-v_{1x}\int_0^x\pi(u_2)dy + \pi(u_2) v_1)_x\right] \end{pmatrix},
	\\\nonumber
J[u, Jv] & = \begin{pmatrix} -\int_0^x \pi\left(u_{2x}\int_0^y \pi(v_2)dy + v_{1xx}u_1\right) dy \\
\left[(u_{1x}\int_0^x\pi(v_2)dy - \pi(v_2) u_1)_x\right] \end{pmatrix},
	\\ \nonumber
- [Ju, Jv] & =  \begin{pmatrix} -\pi(v_2)\int_0^x\pi(u_2)dy +  \pi(u_2) \int_0^x \pi(v_2) dy \\
\left[v_{1xx}\int_0^x \pi(u_2)dy - u_{1xx}\int_0^x\pi(v_2)dy\right]\end{pmatrix},
\end{align}
and the sum of these four equations vanishes after simplification. Thus, $N^J = 0$.
\proofend

\begin{remark}\upshape
Not all the properties listed in Theorem \ref{kahlerth} are independent. For example, suppose that (a) and (c) hold. Then $\omega = g(J \cdot, \cdot)$ satisfies $\nabla \omega = 0$. Provided that the connection $\nabla$ Êis torsion-free, equation (\ref{ddef}) can be rewritten as
$$d\eta(X_0, \dots, X_p) = \sum_{i=0}^p (-1)^i (\nabla_{X_i}\eta)(X_0, \dots, \hat{X}_i, \dots, X_p),$$
and it follows that $\omega$ is closed. On the other hand, suppose that (a), (b), and (f) hold on an almost complex Banach manifold. Then it can be showed that $\nabla J$ must vanish. Indeed, the same argument that is used to prove this in the finite-dimensional setting (see Theorem 11.5 of \cite{M2007}) generalizes to Banach manifolds.
\end{remark}

We can write (\ref{pi2HS}) in the following form suitable for the formulation of weak solutions:
\begin{align}\label{weakpi2HS}
  & \begin{pmatrix} u_t + uu_x \\ \pi(\rho)_t  \end{pmatrix} 
  =  \begin{pmatrix} - \frac{1}{2}A^{-1}\partial_x(u_x^2 + \pi(\rho)^2) \\
-(u\pi(\rho))_x \end{pmatrix}.
\end{align}

\begin{proposition}\label{pigeodesicprop}
Let $s > 5/2$. Let $(\varphi, [\alpha]):J \to K^s$ be a $C^2$-curve where $J \subset \R$ is an open interval and define $(u, \rho)$ by (\ref{urhovarphitft}). Then 
\begin{align}\label{piurhoCC1}
& (u, [\rho]) \in C\left([0, T);  H_0^s(S^1;\R) \times (H^{s-1}(S^1;\R)/\R)\right)
	\\ \nonumber
 & \hspace{3cm}\cap C^1\left([0, T); H_0^{s-1}(S^1;\R) \times (H^{s-2}(S^1;\R)/\R)\right)
\end{align}
and $(\varphi, [\alpha])$ is a geodesic if and only if $(u, [\rho])$ satisfies (\ref{weakpi2HS}).
\end{proposition}
\proofbegin
In terms of $(u, [\rho])$, the geodesic equation 
\begin{align}\label{pigeodesiceq}
  (\varphi_{tt}, [\alpha]_{tt}) = \Gamma_{(\varphi, [\alpha])}((\varphi_t, [\alpha]_t), (\varphi_t, [\alpha]_t))
\end{align}
takes the form
$$\begin{pmatrix} u_t + uu_x \\ [\rho]_t + [u\rho_x] \end{pmatrix} = \Gamma_{(\id, [0])}((u, [\rho]), (u, [\rho])).$$
The first component of this equation is the first component of (\ref{weakpi2HS}).
Applying the map $[\rho] \mapsto \pi(\rho)$ to the second component, we find
$$\pi(\rho)_t + \pi(u\rho_x) = - \pi(u_x\pi(\rho)).$$
Writing this as
$$\pi(\rho)_t = - \pi((u\pi(\rho))_x),$$
we see that it coincides with the second component of (\ref{weakpi2HS}).
\proofend

\section{Complex projective space and the Hopf fibration}\label{hopfsec}\nequation
In this section, we will show that the K\"ahler manifold $K^s$ introduced in Section \ref{kahlersec} is isometric to a subset of complex projective space. Under this isometry, the metric $\langle \cdot, \cdot \rangle$ on $K^s$ is simply the Fubini-Study metric and the K\"ahler structure on $K^s$ corresponds to the canonical K\"ahler structure on $\CP^\infty$. 
Moreover, we will show that the fibration
\begin{align}\label{hopffibration2}
\xymatrix{S_{4\pi}^1 \ar[r] & G^s \ar[d] \\
& K^s }
\end{align}
which arises because the circle $S_{4\pi}^1$ acts on $G^s$ by addition of a constant phase in the second component and the quotient space is $K^s$, corresponds under the above isometries to the Hopf fibration of $S^\infty$ over $\CP^\infty$.

We first define the relevant infinite-dimensional complex projective space in detail. 
Let $\CP^\infty = (L^2(S^1; \C)\setminus\{0\})/\sim$, where $f \sim g$ iff $f(x) = c g(x)$ for some $c \in \C$. For $s>5/2$, let $\CP^{\infty, s}$ denote the elements in $\CP^\infty$ of Sobolev class $H^s$, i.e.
$$\CP^{\infty, s} = (H^s(S^1; \C)\setminus\{0\})/\sim.$$ 
Let $\projV: H^s(S^1; \C)\setminus\{0\} \to \CP^{\infty, s}$ be the natural projection.
We turn $\CP^{\infty, s}$ into a Banach manifold modeled on the space 
$$H_1^s(S^1; \C) := \{f \in H^s(S^1; \C) \; | \; f(0) = 1\}$$ 
as follows: For each $x_0 \in S^1$, we let $W_{x_0}$ be the open subset of $\CP^{\infty, s}$ defined by
$$W_{x_0} = \projV\left(\{f \in \CP^{\infty,s}\; | \; f(x_0) \neq 0 \}\right),$$
and let $\varphi_{x_0}:W_{x_0} \to H_1^s(S^1; \C)$ be the map
$$\varphi_{x_0}: [f] \mapsto \frac{f(\cdot + x_0)}{f(x_0)}.$$
The collection of charts $\{(W_{x_0}, \varphi_{x_0})\}_{x_0 \in S^1}$ covers $\CP^{\infty, s}$. Moreover, for any $x_0, x_1 \in S^1$, the transition map $\varphi_{x_0} \circ \varphi_{x_1}^{-1}$ given by
$$\varphi_{x_0} \circ \varphi_{x_1}^{-1}: f \mapsto \frac{f(\cdot + x_0 - x_1)}{f(x_0 - x_1)}$$
is smooth. This defines the manifold structure on $\CP^{\infty, s}$.

The circle $S_{4\pi}^1$ acts by isometries on the unit sphere $S^{\infty,s} \subset L^2(S^1; \C)$ by multiplication by a constant phase,
$$S_{4\pi}^1 \times S^{\infty,s} \to S^{\infty,s}: f \mapsto fe^{i\alpha_0/2}, \qquad \alpha_0 \in S_{4\pi}^1.$$
The quotient space is $\CP^{\infty, s}$ and the restriction of $\projV$ to $S^{\infty, s}$ is the quotient map for this action.
The action determines, for each $f \in S^{\infty,s}$, an orthogonal splitting of the tangent space $T_fS^{\infty,s} \subset L^2(S^1; \C)$ according to
$$T_fS^{\infty,s} = (T_fS^{\infty,s})^{v} \oplus (T_fS^{\infty,s})^{h},$$
where the vertical and horizontal subspaces are given by
\begin{align}\label{kerTrho}
(T_fS^{\infty,s})^{v} := \ker T_f \projV = \{i c f | c \in \R\},
\end{align}
and
$$(T_fS^{\infty,s})^{h} := (\ker T_f \projV)^\perp = \left\{X \in T_fS^{\infty,s} \;\middle|\; \text{Re} \int if \bar{X} dx = 0\right\},$$
respectively.
The Fubini-Study metric on $\CP^{\infty,s}$ is defined by
$$\langle T\projV(X), T\projV(Y) \rangle_{\text{FS}} = \langle X^{h}, Y^{h} \rangle_{L^2},$$
where $X^h$ and $Y^h$ denote the horizontal components of $X, Y$. The orthogonal projection of $X \in T_gS^{\infty,s}$ onto $(T_gS^{\infty,s})^{h}$ is given by
$$X \mapsto X^{h} = X + ig\text{Im} \int_{S^1} g\bar{X} dx.$$

Let $\mathcal{V}^{s}$ be the image under $\projV$ of the subset $\mathcal{U}^s \subset S^{\infty, s}$ of nowhere vanishing functions defined in (\ref{Usdef}). Let  $\projK: G^s \to K^s$ denote the projection
\begin{align}\label{GtoKprojection}
\projK(\varphi, \alpha) = (\varphi, [\alpha])	
\end{align}
Recall that a smooth submersion $F$ from $M$ to $N$, where $M$Ê andÊ $N$ are (possibly weak) Riemannian manifolds, is a {\it Riemannian submersion} if the restriction of $T_pF$ to the horizontal subspace $(\ker T_pF)^\perp \subset T_pM$ is an isometry onto $T_pN$ for each $p \in M$.

\begin{theorem}\label{hopfth}
Let $s > 5/2$. The map $\Psi:K^s \to \mathcal{V}^{s-1}$ defined by
$$\Psi(\varphi, [\alpha]) \mapsto \left[\sqrt{\varphi_x} e^{i\alpha/2}\right]$$
is a diffeomorphism and an isometry. 
Moreover, the natural projections $\projK: G^s \to K^s$ and $\projV : S^{\infty, s} \to \CP^{\infty, s}$ are Riemannian submersions and the following diagram commutes
\begin{align}\label{commutingdiagram}
\xymatrix{\ar[d]_{\projK} G^s \ar[r]^{\Phi} & \mathcal{U}^{s-1} \ar[d]^{\projV} \\
K^s \ar[r]^{\Psi} & \mathcal{V}^{s-1}}
\end{align}
where $\Phi$ denotes the isometry of Theorem \ref{sphereth}.
\end{theorem}
\proofbegin
The map $\Psi$ is bijective with inverse given by
$$\Psi^{-1}(f) = \left(\int_0^x |f|^2 dx, [2\arg{f(x)}]\right).$$
Since $\Psi$ and $\Psi^{-1}$ are smooth, $\Psi$ is a diffeomorphism. The commutativity of the diagram (\ref{commutingdiagram}) follows by construction.

We next verify that the projection $\projV:S^{\infty, s} \to \CP^{\infty, s}$ is a Riemannian submersion. Smoothness of $\projV$ can be verified in local charts. For example, the local representative of $\projV$ with respect to the charts determined by $\varphi_{x_0}$ and the stereographic projection $\sigma_S$ (cf. (\ref{sSdef})) is the smooth map
$$\varphi_{x_0} \circ \projV \circ \sigma_S^{-1}: h \mapsto \frac{2h(\cdot + x_0) - \|h\|_{L^2}^2 + 1}{2h(x_0) - \|h\|_{L^2}^2 + 1}.$$
By definition of the Fubini-Study metric, $\projV$ is a Riemannian submersion.

The projection $\projK: G^s \to K^s$ is also a Riemannian submersion. Indeed, smoothness of $\projK$ is immediate, and for each $(\varphi, \alpha) \in G^s$, $\projK$ determines the splitting
$$T_{(\varphi, \alpha)}G = (T_{(\varphi, \alpha)}G)^{v} \oplus (T_{(\varphi, \alpha)}G)^{h},$$
where the vertical and horizontal subspaces are defined by
$$(T_{(\varphi, \alpha)}G)^{v} := \ker{T_{(\varphi, \alpha)}\projK} = \{(0, U_2)\; |\; U_2 \text{ is a constant function}\},$$
and
$$(T_{(\varphi, \alpha)}G)^{h}  := (\ker{T_{(\varphi, \alpha)}\projK})^\perp = \{(U_1, U_2)\; |\; \pi(U_2) = U_2\},$$
respectively. The orthogonal projections onto the vertical and horizontal subspaces are given by
$$(U_1, U_2) \mapsto (U_1, U_2)^{v} = \left(0, \int_{S^1} U_2 \varphi_x dx\right)$$
and
\begin{align}\label{horliftU}
  (U_1, U_2) \mapsto (U_1, U_2)^{h} = (U_1, \pi(U_2)),
\end{align}
respectively. Let $U^{h} = (U_1, U_2)$ and $V^{h} = (V_1, V_2)$ be horizontal vectors in $T_{(\varphi, \alpha)}G^s$. Then, since $T\projK(U_1, U_2) = (U_1, [U_2])$,
\begin{align*}
\langle U^{h}, V^{h} \rangle_{(\varphi, \alpha)}
&= \frac{1}{4} \int_{S^1} \left(\frac{U_{1x}V_{1x}}{\varphi_x} + \pi(U_2)\pi(V_2)\varphi_x\right)dx
	\\
&= \langle (U_1, [U_2]), (V_1, [V_2]) \rangle_{(\varphi, [\alpha])}
	\\
&=  \langle T\projK(U^{h}), T\projK(V^{h}) \rangle_{(\varphi, [\alpha])},
\end{align*}
showing that $\projK$ is a Riemannian submersion.

Since both $\projV$ and $\projK$ are Riemannian submersions and $\Phi$ is an isometry, it follows from the commuting diagram (\ref{commutingdiagram}) that $\Psi$ also is an isometry.
\proofend

\begin{corollary}
Let $e = (\id, [0])$ denote the identity element in $K^s$. The curvature tensor $R$ on $K^s$ satisfies
\begin{align}\label{Kcurvature}
\langle R(u, v)v, u\rangle_{e}
= \|u\|_{e}^2\|v\|_{e}^2 - \langle u, v\rangle_{e}^2 + 3\omega(u,v)^2,
\end{align}
where $u = (u_1, [u_2])$ and $v = (v_1, [v_2])$ are elements in $T_{e}K^s$. In particular, the sectional curvature 
$$\sect(u,v) = \frac{\langle R(u, v)v, u\rangle_{e}}{\|u\|_e^2\|v\|_e^2 - \langle u, v\rangle_e^2}$$
satisfies
\begin{align}\label{1sec4}
1 \leq \sect(u,v) \leq 4, \qquad u,v \in T_{e}K^s,
\end{align}
and $\sect(u,v) = 4$ if and only if $Ju$ is a multiple of $v$.
\end{corollary}
\proofbegin
We will give three different proofs of (\ref{Kcurvature}). The first proof is the most elegant and relies on Theorem \ref{hopfth}. According to this theorem, $\projK:G^s \to K^s$ is a Riemannian submersion. The O'Neill formula for Banach manifolds (see \cite{Lang}) implies that
$$\langle R(X,Y)Y, X \rangle_{K^s} = \langle R_{G^s}(X^h, Y^h)Y^h, X^h \rangle_{G^s} + \frac{3}{4}\bigl\|\bigl[X^h, Y^h\bigr]^{v}\bigr\|_{G^s}^2,$$
where $X^h, Y^h$ denote the horizontal lifts of two orthonormal vector fields $X, Y$ on $K^s$ and $R_{G^s}$ denotes the curvature tensor on $G^s$.
In view of Corollary \ref{curvcor} and equation (\ref{horliftU}) this yields
$$\langle R(u,v)v, u \rangle_{K^s} = \|u\|_e^2\|v\|_e^2 - \langle u, v\rangle_e^2  + \frac{3}{4} \bigl\|[(u_1, \pi(u_2)), (v_1, \pi(v_2))]^{v}\bigr\|_{G^s}^2.$$
Since 
$$[(u_1, \pi(u_2)), (v_1, \pi(v_2))]^{v} 
= \begin{pmatrix} v_{1x}u_1 - u_{1x}v_1 \\ v_{2x}u_1 - u_{2x}v_1 \end{pmatrix}^{v}
= \begin{pmatrix} 0 \\ \int_{S^1}(v_{2x}u_1 - u_{2x}v_1)dx\end{pmatrix},$$
we find (\ref{Kcurvature}).

The second proof utilizes the following formula for the curvature of a Lie group with a right-invariant metric, see \cite{Abook}:
\begin{align}\label{arnoldcurvformula}
  \langle R(u,v)v, u \rangle_e = \langle \delta, \delta \rangle_e + \langle [u,v], \beta \rangle_e - \frac{3}{4}\langle [u,v], [u,v]\rangle_e - \langle B(u,u), B(v,v) \rangle_e,
\end{align}
where 
$$\delta := \frac{1}{2} \left(B(u,v) + B(v,u) \right), \qquad \beta := \frac{1}{2} \left(B(u,v) - B(v,u) \right),$$
and the bilinear map $B$ is defined by
$$\langle B(u,v), w \rangle_e = \langle u, [v, w]\rangle_e.$$ 
Using expression (\ref{Ksbracket}) for the Lie bracket, we infer that 
$$B(u,v) = \begin{pmatrix} A^{-1}(u_{1xx}v_{1x} + (u_{1xx}v_1)_x - \pi(u_2)v_{2x}) \\
- [(\pi(u_2)v_1)_x] \end{pmatrix}$$
and tedious computations show that (\ref{arnoldcurvformula}) reduces to (\ref{Kcurvature}).

The third proof employs the local expression for the curvature tensor $R$ in terms of the Christoffel map:
\begin{align*}
R_p(U, V)W = & D_1\Gamma_{p}(W, U)V - D_1\Gamma_{p}(W, V)U
  		 \\
 &  + \Gamma_{p}(\Gamma_{p}(W, V), U) - \Gamma_{p}(\Gamma_{p}(W, U), V)
\end{align*}
where $D_1$ denotes differentiation with respect to $p$.
Proceeding as in the proof of Proposition 5.1 of \cite{EKL2011}, it can be shown that
\begin{align} \nonumber
\langle R(u, v)v, u\rangle_{e} =&\;\langle\Gamma(u,v), \Gamma(u,v) \rangle_e -\langle\Gamma(u,u), \Gamma(v,v) \rangle_e
	\\ \label{Suv2}
&- \left\langle\begin{pmatrix}
  u_{1x}v_1 \\
  u_{2x}v_1 
\end{pmatrix}, \Gamma(u,v) \right\rangle_e
+\left\langle
\begin{pmatrix}
  u_{1x}u_1 \\
  u_{2x}u_1 
\end{pmatrix}, \Gamma(v,v) \right\rangle_e
	\\ \nonumber
&+\langle-\Gamma(v_xv_1,u)-\Gamma(v,u_xv_1)+2\Gamma(v_xu_1,v), u\rangle_e.
\end{align}
Letting $\mu(f) = \int_{S^1} f(x) dx$ denote the mean of a function $f:S^1 \to \R$, long but straightforward computations yield
\begin{align*}
& \langle \Gamma(u,v), \Gamma(u,v) \rangle_e - \langle \Gamma(u,u), \Gamma(v,v) \rangle_e
= \|u\|^2\|v\|^2 - \langle u, v\rangle_e^2
	\\
& - \frac{1}{16}\left(\mu(u_1v_{2x})^2 + \mu(v_1u_{2x})^2\right) - \frac{1}{8}\mu(v_1u_{2x})\mu(u_1v_{2x})
+ \frac{1}{4}\mu(u_1u_{2x})\mu(v_1v_{2x})
\end{align*}
and
\begin{align*} \nonumber
&-\left\langle\begin{pmatrix}
  u_{1x}v_1 \\
  u_{2x}v_1 
\end{pmatrix}, \Gamma(u,v)\right\rangle_e
+\left\langle
\begin{pmatrix}
  u_{1x}u_1 \\
  u_{2x}u_1 
\end{pmatrix}, \Gamma(v,v) \right\rangle_e
	\\ \nonumber
&+\langle-\Gamma(v_xv_1,u)-\Gamma(v,u_xv_1)+2\Gamma(v_xu_1,v), u\rangle_e
	\\
& \hspace{2cm} = \frac{1}{4}\left\{\mu(u_1v_{2x})^2 + \mu(v_1u_{2x})^2 - \mu(v_1v_{2x}) \mu(u_1u_{2x}) - \mu(v_1u_{2x})\mu(u_1v_{2x})\right\},
\end{align*}
which again leads to formula (\ref{Kcurvature}).

To prove (\ref{1sec4}) we assume that $u$ and $v$ are orthogonal. Then
$$\frac{3\omega(u,v)^2}{\|u\|_e^2\|v\|_e^2 - \langle u, v\rangle_e^2}
 = \frac{3g(Ju,v)^2}{\|u\|_e^2\|v\|_e^2} \leq \frac{3\|Ju\|_e^2\|v\|_e^2}{\|u\|_e^2\|v\|_e^2} = 3,$$ 
with equality iff $Ju$ is a multiple of $v$.
\proofend

\section{Conclusions and remarks}\label{conclusionsec}\nequation
As was noted in the introduction, the flow of the classical Kepler problem in celestial mechanics is equivalent to the geodesic flow on a sphere. In this paper, we showed that the flow of the Hunter-Saxton system (\ref{2HS}) is also equivalent to the geodesic flow on a sphere, namely, to the geodesic flow on (a subset of) the unit sphere in $L^2(S^1; \C)$. 
Utilizing this geometric picture, we were able to integrate equation (\ref{2HS}) explicitly. Moreover, an infinite-dimensional example of a Hopf fibration was obtained by considering the restriction of (\ref{2HS}) to solutions $(u,\rho)$ where $\rho$ has zero mean. The restricted equation describes the geodesic flow on an infinite-dimensional complex projective space with a natural K\"ahler structure.

The Hunter-Saxton equations (\ref{HS}) and (\ref{2HS}) are the special cases when $M = S^1$ of two more general equations introduced in \cite{KMLP2011} and \cite{LY2011}, which are defined for any compact Riemannian manifold $M$. The first of these equations describes the geodesic flow on the unit sphere in $L^2(M; \R)$ \cite{KMLP2011}, whereas no such sphere interpretation is as yet known for the second equation (although it is known that the underlying space has constant positive curvature \cite{LY2011}). The considerations of this paper suggest that the flow of the equation in \cite{LY2011} is equivalent to the geodesic flow on (a subset of) the unit sphere in $L^2(M; \C)$. The details of this construction will be considered elsewhere.

We emphasize that the fact that the underlying spaces for the above four geodesic equations (the equation in \cite{LY2011} and its three special cases (\ref{HS}), (\ref{2HS}), and the equation in \cite{KMLP2011}) have constant positive curvature is exceptional in the context of PDEs that arise as geodesic equations. In most cases, the curvature takes on both signs (one exception is the inviscid Burgers equation in one space dimension, for which the curvature is everywhere non-negative; in higher dimensions this is no longer true).


\bigskip
\noindent
{\bf Acknowledgement} {\it The author acknowledges support from the EPSRC, UK.}

\bibliographystyle{plain}

\end{document}